\documentclass[11pt,reqno]{amsart}

\usepackage{amsmath}
\usepackage{amssymb}
\usepackage{amsthm}

\newtheorem{theorem}{Theorem}[section]
\newtheorem{prop}[theorem]{Proposition}

\theoremstyle{definition}
\newtheorem{defn}[theorem]{Definition}
\theoremstyle{remark}
\newtheorem*{rem}{Remark}
\numberwithin{equation}{section}

\newcommand{\Gg}{\mathfrak{g}}    
\newcommand{\Gh}{\mathfrak{h}}

\newcommand{\Gz}{\mathfrak{z}}

\begin{document}

\title[Quantizations of some Poisson--Lie groups]
{Quantizations of some Poisson--Lie groups: The bicrossed product construction}

\author{Byung-Jay Kahng}
\date{}
\address{Department of Mathematics and Statistics\\ University of Nevada\\
Reno, NV 89557}
\email{bjkahng@unr.edu}

\begin{abstract}
By working with several specific Poisson--Lie groups arising from Heisenberg
Lie bialgebras and by carrying out their quantizations, a case is made for
a useful but simple method of constructing locally compact quantum groups.
The strategy is to analyze and collect enough information from a Poisson--Lie
group, and using it to carry out a ``cocycle bicrossed product construction''.
Constructions are done using multiplicative unitary operators, obtaining
$C^*$-algebraic, locally compact quantum (semi-)groups.
\end{abstract}
\maketitle

{\sc Introduction.}
Typically, the most commonly employed method of constructing specific
examples of quantum groups is the method of ``generators and relations''.
This is certainly the case in the purely algebraic setting of quantized
universal enveloping (QUE) algebras.  Even in the $C^*$-algebra setting,
compact quantum groups are usually constructed in this way (See, for
instance, \cite{Wr1}.).

However, when one wishes to construct a non-compact quantum group,
this method is not so useful: The generators (essentially the
coordinate functions of a group) tend to be unbounded, which gives
rise to various technical difficulties.  There does exist ways to
handle the difficulties (See \cite{Wr4}, where Woronowicz works with
the notion of unbounded operators ``affiliated'' with $C^*$-algebras),
but in general, it is usually better to look for some other methods
of construction.

One useful approach not relying on the generators is via the method
of deformation quantization.  Here, the aim is to deform the
(commutative) algebra of functions on a Poisson manifold, in
the direction of the Poisson bracket (See \cite{BFFL}, \cite{Vy}.
For the ``strict deformation quantization'' in the $C^*$-algebra
framework, see also \cite{Rf1}, \cite{Rf4}.).  One should note,
however, that this is just a ``spatial'' deformation, in the sense
that the deformation is carried out only for the algebra structure.
To obtain a quantum group, one begins with a suitable Poisson--Lie
group $G$ (a Lie group equipped with a compatible Poisson bracket)
and perform the deformation quantization on the function space
$C_0(G)$---for both its algebra and coalgebra structures.

Some of the non-compact quantum groups obtained by deformation
quantization are \cite{Rf5}, \cite{SZ}, \cite{Rf7}, \cite{Zk2},
\cite{BJKp2}.  In these examples, the information at the level
of Poisson--Lie groups or Lie bialgebras plays a key role in
constructing the quantum groups and their structures.  Naturally,
there exists a very close relationship between a quantum group
obtained in this way and its Poisson--Lie group counterpart.

This point is quite helpful in working with the quantum group.
For instance, as for the example considered by the author in
\cite{BJKp2}, the information from the classical (Poisson) level
was useful not only in the construction of the quantum group
but also in studying its representation theory, in relation
to the dressing orbits \cite{BJKhj}, \cite{BJKppdress}.

On the other hand, despite many advantages, there are some
drawbacks to the method of deformation quantization, especially
when one wishes to carry it out in the $C^*$-algebra setting:
Jumping from the classical level of Poisson--Lie groups to the
$C^*$-algebraic quantum group level is not necessarily an easy
task.  Even with the guides provided from Poisson level, the
actual construction of the structure maps like comultiplication,
antipode, or Haar weight often should be done using different
means.  Among the useful tools is the notion of multiplicative
unitary operators (in the sense of Baaj and Skandalis \cite{BS}).

Considering the drawbacks to the geometric approach above, we
turn to a more algebraic method of constructing locally compact
quantum groups, via the framework of (cocycle) bicrossed products.
This goes back to the problem of group extensions in the Kac
algebra setting (see \cite{ES} for a survey on Kac algebras),
and was made systematic by Majid \cite{Mj1}, \cite{Mj0}.  Here,
one begins with a certain ``matched pair'' of groups (or more
generally, locally compact quantum groups) and build a larger
quantum group as a bicrossed product, possibly with a cocycle.
Baaj and Skandalis has a version of this in Section 8 of \cite
{BS}.  For a comprehensive treatment about this framework, see
\cite{VV}.

The best aspect about the bicrossed product method is that
it is relatively simple, while sufficiently general to
include many special cases.  However, as is the case for
any general method, having the framework is not enough to
construct actual and specific examples: One needs to have
a specific matched pair, together with a compatible cocycle,
for this method to work.

So we propose here to combine the advantages of the
``geometric'' (deformation quantization) method and the
``algebraic'' framework of cocycle bicrossed products.
That is, we first begin with a Poisson--Lie group and
analyze its Poisson structure.  The Poisson data will
help us obtain a suitable matched pair and a compatible
cocycle.  Then we perform the cocycle bicrossed product
construction.

This program is usually best for constructing solvable-type
quantum groups.  It is because crossed products often model
quantized spaces (For instance, the ``Weyl algebra'',
$C_0(\mathbb{R}^n)\rtimes_{\tau}\mathbb {R}^n$ with $\tau$
being the translation, is the quantized phase space \cite
{Fo}.).  With some adjustments, it can be also useful to
construct even wider class of quantum groups.  Moreover,
having a close connection with the Poisson--Lie group
enables us to take advantage of its geometric data in
further studying the quantum group as well as in applications.

Our plan in this paper is to illustrate this program through
examples, using Poisson--Lie groups associated with several
Heisenberg-type Lie bialgebras.  Three specific cases are
considered: In Case (2), we re-construct our earlier example
from \cite{BJKp2}, which quantizes a certain non-linear
Poisson bracket.  Case (1) is about examples from \cite{Rf5},
\cite{SZ}, \cite{VD}, which are related with a linear Poisson
bracket (so a little simpler).  Case (3) is similar to the
example given in \cite{EV}, but is more general.  Afterwards,
we give more constructions of similar-flavored examples.

To keep the presentation coherent and simple, we will not stray
too much away from the Heisenberg-type Lie bialgebras and their
quantum counterparts.  In this article, the focus is not on
giving genuinely new examples, but on establishing a simple
but quite useful method of constructing specific quantum groups.
In this way, we make a case that the geometric, deformation
quantization method and the algebraic, bicrossed product method
are very much compatible.  Examples constructed with the same
program but coming from different Poisson--Lie groups will be
presented in our future work.

The paper is organized as follows.  In Section 1, we discuss
some specific Poisson structures coming from Heisenberg Lie
bialgebras.  Three specific cases are considered.  In Section 2,
we carry out the quantizations of the cases from Section 1.
By analyzing the Poisson brackets, we obtain, for each of the
three cases, a matched pair and a cocycle.  These data help
us to construct a multiplicative unitary operator, which
represents the cocycle bicrossed product construction.

More examples are given in Section 3, by slightly modifying
the results obtained in Section 2.  Along the way, we will
make frequent comparisons between the examples in Sections
2 and 3, and several other examples obtained elsewhere using
different methods.

Section 4 is Appendix, showing that the Poisson structures
considered in Section 1 actually arise from certain ``classical
$r$-matrices''.  Remembering that a ``quantum $R$-matrix'' type
operator played a significant role in the representation theory
of our earlier example (\cite{BJKp2}, \cite{BJKhj}), this is
a useful knowledge.

\section{Lie bialgebra structures on a Heisenberg Lie algebra.
The Poisson--Lie groups}

Let $H$ be the $(2n+1)$-dimensional Heisenberg Lie group.  Its
underlying space is $\mathbb{R}^{2n+1}$ and the multiplication
on it is given by
$$
(x,y,z)(x',y',z')=\bigl(x+x',y+y',z+z'+\beta(x,y')\bigr),
$$
for $x,x',y,y'\in\mathbb{R}^n$ and $z,z'\in\mathbb{R}$.  Here,
$\beta(\ ,\ )$ denotes the ordinary inner product.

Its Lie algebra counterpart is the Heisenberg Lie algebra $\Gh$.
It is generated by $\mathbf{x_i},\mathbf{y_i} (i=1,\dots,n),
\mathbf{z}$, with the following relations:
$$
[\mathbf{x_i},\mathbf{y_j}]=\delta_{ij}\mathbf{z}, 
\quad [\mathbf{z},\mathbf{x_i}]=[\mathbf{z},\mathbf{y_i}]=0.
$$

For convenience, we will identify $H\cong\Gh$ as vector spaces.
This is possible since $H$ is an exponential solvable Lie group
(it is actually nilpotent).  We will understand that $x=x_1
\mathbf{x_1}+\dots+x_n\mathbf{x_n}$, and similarly for the other
variables.  And, we choose a Lebesgue measure on $H\cong\Gh$,
which is a Haar measure for $H$.

For a Heisenberg Lie group, all the possible compatible Poisson 
brackets on it have been classified by Szymczak and Zakrzewski
\cite{SZ}.  Among these, we will specifically look at the
following, simpler cases (See Definition \ref{PoissonH} below.).
The Poisson brackets are described in terms of the cobrackets
$\delta:\Gh\to\Gh\wedge\Gh$, which are 1-cocycles with respect
to the adjoint representation.  It is known from general theory
that specifying in such a way a ``Lie bialgebra'' structure,
$(\Gh,\delta)$, is equivalent to giving an explicit formula for
the Poisson bracket on $H$.  See \cite{LW}.

\begin{defn}\label{PoissonH}
\begin{enumerate}
\item Consider $\delta_1:\Gh\to\Gh\wedge\Gh$ defined by
$$
\delta_1(\mathbf{x_j})=\lambda\mathbf{x_j}\wedge\mathbf{z},\quad
\delta_1(\mathbf{y_j})=-\lambda\mathbf{y_j}\wedge\mathbf{z},\quad
\delta_1(\mathbf{z})=0. 
$$
Here $\lambda\in\mathbb{R}$.  To obtain a nontrivial map, we let 
$\lambda\ne0$.
\item Let $\lambda\ne0$ again, and let $\delta_2:\Gh\to\Gh\wedge\Gh$
be defined by
$$
\delta_2(\mathbf{x_j})=\lambda\mathbf{x_j}\wedge\mathbf{z},\quad
\delta_2(\mathbf{y_j})=\lambda\mathbf{y_j}\wedge\mathbf{z},\quad
\delta_2(\mathbf{z})=0. 
$$
\item Let $(J_{ij})$ be a skew, $n\times n$ matrix ($n\ge2$), and
let $\delta_3:\Gh\to\Gh\wedge\Gh$ be defined by
$$
\delta_3(\mathbf{x_j})=0,\quad
\delta_3(\mathbf{y_j})=\sum_{i=1}^n J_{ij}\mathbf{x_i}\wedge
\mathbf{z},\quad
\delta_3(\mathbf{z})=0. 
$$
\end{enumerate}
\end{defn}

We do not plan to give here an explicit proof that these are indeed
Lie bialgebra structures on $\Gh$ giving us the compatible Poisson
brackets on $H$.  Instead, we can refer to Theorem 2.2 of \cite{SZ},
and in the case of $\delta_2$ above, a careful discussion was given
in Section 1 of \cite{BJKp2}.  See also Section 4 (Appendix), where
we show that they arise from certain ``classical $r$-matrices''.

Corresponding to each of these Poisson brackets, we can define a Lie
bracket on the dual space ${\Gh}^*$ of $\Gh$ by $[\ ,\ ]=\delta^*:
{\Gh}^*\wedge{\Gh}^*\to{\Gh}^*$.  That is, $[\mu,\nu]$ is defined by
\begin{equation}\label{(delta^*)}
\bigl\langle[\mu,\nu],X\bigr\rangle=\bigl\langle\delta^*(\mu\otimes\nu),
X\bigr\rangle=\bigl\langle\mu\otimes\nu,\delta(X)\bigr\rangle,
\end{equation}
where $X\in\Gh$, $\mu,\nu\in{\Gh}^*$, and $\langle\ ,\ \rangle$ is
the dual pairing between ${\Gh}^*$ and $\Gh$.  In this way, we obtain
the following ``dual'' Lie algebra for each of the cases.  The proof
is straightforward.

\begin{prop}\label{dualLiealgebra}
Let $\Gg={\Gh}^*$ be spanned by $\mathbf{p_i},\mathbf{q_i}
(i=1,\dots,n),\mathbf{r}$, which form the dual basis of $\mathbf
{x_i},\mathbf{y_i} (i=1,\dots,n),\mathbf{z}$.
\begin{enumerate}
\item On $\Gg$, define the Lie algebra relations for the basis
vectors as follows:
$$
[\mathbf{p_i},\mathbf{q_j}]=0,  \quad
[\mathbf{p_i},\mathbf{r}]=\lambda\mathbf{p_i}, \quad
[\mathbf{q_i},\mathbf{r}]=-\lambda\mathbf{q_i}.  
$$
Then $\Gg$ is the Poisson dual of the Lie bialgebra $(\Gh,\delta_1)$.
\item On $\Gg$, define the Lie algebra relations for the basis vectors
by
$$
[\mathbf{p_i},\mathbf{q_j}]=0, \quad
[\mathbf{p_i},\mathbf{r}]=\lambda\mathbf{p_i}, \quad
[\mathbf{q_i},\mathbf{r}]=\lambda\mathbf{q_i}.
$$
This is the Poisson dual of the Lie bialgebra $(\Gh,\delta_2)$.
\item On $\Gg$, define the Lie algebra relations by
$$
[\mathbf{p_i},\mathbf{q_j}]=0, \quad
[\mathbf{p_i},\mathbf{r}]=\sum_{j=1}^n J_{ij}\mathbf{q_j},  \quad
[\mathbf{q_i},\mathbf{r}]=0.   \notag
$$
This is the Poisson dual of the Lie bialgebra $(\Gh,\delta_3)$.
\end{enumerate}
\end{prop}

Each of the dual Lie algebras $\Gg$ is actually a Lie bialgebra,
whose cobracket $\theta:\Gg\to\Gg\wedge\Gg$ is the dual map of
the Lie bracket on $\Gh$.  This situation is exactly same as
in equation \eqref{(delta^*)}.  In other words, $\theta:\Gg\to
\Gg\wedge\Gg$ is defined by its values on the basis vectors of
$\Gg$ as follows:
$$
\theta(\mathbf{p_i})=0,\quad\theta(\mathbf{q_i})=0,
\quad\theta(\mathbf{r})=\sum_{i=1}^n(\mathbf{p_i}\otimes
\mathbf{q_i}-\mathbf{q_i}\otimes\mathbf{p_i})=\sum_{i=1}^n
(\mathbf{p_i}\wedge\mathbf{q_i}). 
$$

We thus have the (Poisson dual) Lie bialgebra $(\Gg,\theta)$,
for each of the Heisenberg Lie bialgebras in Definition \ref
{PoissonH}.  Let us now consider the corresponding Poisson--Lie
groups $G$ (dual to the Heisenberg Lie group), together with
their Poisson brackets.  As before, we will understand $p=p_1
\mathbf{p_1}+p_2\mathbf{p_2}+\cdots+p_n\mathbf{p_n}$, and
similarly for the other variables (This explains the notation
we use in (3) below.).

\begin{prop}\label{dualLiegroups}
\begin{enumerate}
\item Let $G$ be the $(2n+1)$-dimensional Lie group, whose
underlying space is $\mathbb{R}^{2n+1}$ and the multiplication
law is defined by
$$
(p,q,r)(p',q',r')=(e^{\lambda r'}p+p',e^{-\lambda r'}q+q',r+r').
$$
It is the Lie group corresponding to $\Gg$ from Proposition
\ref{dualLiealgebra} (1).  The Poisson bracket on $G$ is given
by the expression
$$
\{\phi,\psi\}(p,q,r)=r\bigl(\beta(x,y')-\beta(x',y)\bigr),
$$
for $\phi,\psi\in C^{\infty}(G)$.  Here $d\phi(p,q,r)=(x,y,z)$
and $d\psi(p,q,r)=(x',y',z')$, which are naturally considered
as elements of $\Gh$.
\item Let $G$ be the $(2n+1)$-dimensional Lie group, whose
underlying space is $\mathbb{R}^{2n+1}$ and together with the
multiplication law
$$
(p,q,r)(p',q',r')=(e^{\lambda r'}p+p',e^{\lambda r'}q+q',r+r').
$$
It is the Lie group corresponding to $\Gg$ from Proposition
\ref{dualLiealgebra} (2).  The Poisson bracket on $G$ is given
by
$$
\{\phi,\psi\}(p,q,r)=\left(\frac{e^{2\lambda r}-1}{2\lambda}
\right)\bigl(\beta(x,y')-\beta(x',y)\bigr),
$$
for $\phi,\psi\in C^{\infty}(G)$.  Here again, we use the
natural identification of $d\phi(p,q,r)=(x,y,z)$ and $d\psi
(p,q,r)=(x',y',z')$ as elements of $\Gh$.
\item Let $G$ be the $(2n+1)$-dimensional Lie group, together
with the multiplication law
$$
(p,q,r)(p',q',r')=\left(p+p',q+q'+r'\sum_{i,j=1}^n J_{ij}p_i
\mathbf{q_j},r+r'\right).
$$
For it to be non-trivial, we need $n\ge2$.  This gives us
the Lie group corresponding to $\Gg$ from Proposition \ref
{dualLiealgebra} (3).  The Poisson bracket on $G$ is given
by
$$
\{\phi,\psi\}(p,q,r)=r\bigl(\beta(x,y')-\beta(x',y)\bigr)
+\frac{r^2}{2}\sum_{k,j=1}^n J_{kj}(y_jy'_k-y_ky'_j),
$$
for $\phi,\psi\in C^{\infty}(G)$.  Again, $d\phi(p,q,r)
=(x,y,z)$ and $d\psi(p,q,r)=(x',y',z')$, viewed as elements
of $\Gh$.
\end{enumerate}
\end{prop}

\begin{proof}
Constructing $G$ from $\Gg$ is rather straightforward.  In each
of the three cases, $G$ is a (connected and simply connected)
exponential solvable Lie group corresponding to $\Gg$.  As before,
we can identify $G\cong\Gg$ as vector spaces.  Note that the
definitions of the group multiplications are chosen in such a
way that an ordinary Lebesgue measure becomes the Haar measure
for $G$ (in particular, for Case (2)).

To find the expression for the Poisson bracket, we follow the
standard procedure:  First, consider $\operatorname{Ad}:G\to
\operatorname{Aut}(\Gg)$, the adjoint representation of $G$ on
$\Gg$.  We then look for a map $F:G\to\Gg\wedge\Gg$, that is a
group 1-cocycle on $G$ for the $\operatorname{Ad}$ representation
and whose derivative at the identity element, $dF_e$, coincides
with $\theta$ above.  Note that since $\theta$ depends only on
the $r$-variable, so should $F$.  In other words, we look for
a map $F$ such that
\begin{align}
F(r+r')&=F(r)+\operatorname{Ad}_{(0,0,r)}\bigl(F(r')\bigr),
\notag \\
dF_{(0,0,0)}(r)&=\theta(r)=r\sum_{k=1}^n(\mathbf{p_k}\wedge
\mathbf{q_k}).  \notag
\end{align}
Once we have the 1-cocycle $F$, the Poisson bivector field is
then obtained by the right translation of $F$.

It is true that integrating $\theta$ to $F$ is not always easy.
However, it is not too difficult in our three cases above, due to
our Lie bialgebra structures being rather simple.  In particular,
for Case (2), the computation was given in the proof of Theorem 2.2
in \cite{BJKp2}.  Case (1) is similar but easier, since the map
$F$ (and the Poisson bracket) is linear.

As for Case (3), note first that the representation $\operatorname
{Ad}$ sends the basis vectors of $\Gg$ as follows:
\begin{align}
&\operatorname{Ad}_{(0,0,r)}(\mathbf{p_k})=(0,0,r)\mathbf{p_k}(0,0,-r)
=\mathbf{p_k}-r\sum_{j=1}^n J_{kj}\mathbf{q_j},  \notag \\
&\operatorname{Ad}_{(0,0,r)}(\mathbf{q_k})=\mathbf{q_k},\qquad
\operatorname{Ad}_{(0,0,r)}(\mathbf{r})=\mathbf{r}.  \notag
\end{align}
Considering the requirements for the map $F$ given above, we obtain
the following expression for $F$:
$$
F(p,q,r)=F(r)=r\sum_{k=1}^n(\mathbf{p_k}\wedge\mathbf{q_k})-\frac
{r^2}{2}\sum_{k,j=1}^n (J_{kj}\mathbf{q_j}\wedge\mathbf{q_k}).
$$
Since the right translations are ${R_{(p,q,r)}}_*(\mathbf{p_k})
=\mathbf{p_k}+r\sum_{j=1}^n J_{kj}\mathbf{q_j}$ and since
${R_{(p,q,r)}}_*(\mathbf{q_k})=\mathbf{q_k}$, we thus have
the expression for our Poisson bracket:
\begin{align}
\{\phi,\psi\}(p,q,r)&=\bigl\langle {R_{(p,q,r)}}_* F(p,q,r),
d\phi(p,q,r)\wedge d\psi(p,q,r)\bigr\rangle  \notag \\
&=r\bigl(\beta(x,y')-\beta(x',y)\bigr)+\frac{r^2}{2}
\sum_{k,j=1}^n J_{kj}(y_jy'_k-y_ky'_j).  \notag
\end{align}
\end{proof}

\begin{rem}
Cases (1) and (2) look almost the same, and the difference may
look rather innocent.  However, Case (1) gives us a linear Poisson
bracket, while Case (2) is non-linear.  Note also that the group
$G$ is unimodular in Case (1), while $G$ is non-unimodular in
Case (2).  It turns out that Case (2) is technically deeper, while
having richer properties: For instance, in working with the Haar
weight and in representation theory of its quantum group counterpart
(See \cite{BJKp2}, \cite{BJKppha}, \cite{BJKhj}, \cite{BJKppdress}.).
Meanwhile, Case (3) gives another non-linear Poisson bracket (in
this situation, $G$ is unimodular).  
\end{rem}

\section{Construction of quantum groups}

Now that we have described our Poisson--Lie groups, let us
construct their quantum group counterparts.  But first, we
should mention that these cases are not totally new, having
been studied elsewhere previously (though Case~(3) will be
somewhat new).  As we noted in Introduction, our real focus
is on illustrating our improved approach of using the Poisson
data to obtain quantum groups, via ``cocycle bicrossed products''.

Considering this, it will be sufficient to just describe
appropriate multiplicative unitary operators.  As is known
in general theory (\cite{BS}, \cite{Wr7}, \cite{KuVa}),
having a ``regular'' (or more generally, ``manageable'')
multiplicative unitary operator gives rise to a
$C^*$-bialgebra, which is really a quantum semigroup.

As for the theory on cocycle bicrossed products in the
$C^*$-algebraic quantum group setting, we refer \cite{VV}.
On the other hand, since we are planning to work with
multiplicative unitary operators, our approach will be
actually closer to that given in Section 8 of \cite{BS}.

In the below, we treat separately the three cases we
described in the previous section.  Using the Poisson
geometric data as a guide, we will find a suitable matched
pair and a compatible cocycle.  We will use these informations
to construct a multiplicative unitary operator, giving rise
to a $C^*$-bialgebra having the structure of a cocycle
bicrossed product.

\subsection{Cases (1) and (2)}
Finding the multiplicative unitary operators for Cases (1) 
and (2) goes essentially the same way.  Since Case (2) is
the more complicated one between the two, we will look at
this case more carefully.  Compare with the construction
procedure given in Sections 2 and 3 of \cite{BJKp2}, where
the approach relied much more on Poisson geometry and the
deformation process.

Since we will use the non-linear expression $(e^{2\lambda r}-1)/
2\lambda$ quite often, let us give it the special notation,
$\eta_{\lambda}(r)$.  Note that if $\lambda=0$, it degenerates
into the linear expression $\eta_{\lambda=0}(r)=r$.

From now on, let the group $G$ and the Poisson bracket on it
be as described in Proposition \ref{dualLiegroups} (2).  It is
the dual Poisson--Lie group of $H$, corresponding to the Lie
bialgebra $(\Gh,\delta_2)$.  Meanwhile, let $Z=\bigl\{(0,0,z):
z\in\mathbb{R}\bigr\}$ be the center of $H$.  Its Lie algebra
counterpart is denoted by $\Gz(\subseteq\Gh)$.

From the expression of the Poisson bracket, we obtain the
following continuous field of group cocycles for $H/Z$:

\begin{prop}\label{cocycle}
Let $r\in{\Gh^*}/{\Gz^{\bot}}$.  Define the map $\sigma^r:H/Z
\times H/Z\to\mathbb{T}$ by
$$
\sigma^r\bigl((x,y),(x',y')\bigr)=\bar{e}\bigl[\eta_{\lambda}(r)
\beta(x,y')\bigr],
$$
where $e(t)=e^{2\pi it}$, so $\bar{e}(t)=e^{-2\pi it}$.  Then
each $\sigma^r$ is a $\mathbb{T}$-valued, normalized group
cocycle for $H/Z$.  Moreover, $r\mapsto\sigma^r$ forms a
continuous field of cocycles.
\end{prop}

\begin{rem}
Verifying the cocycle identity is straightforward, and the
continuity is also clear.  The point is that our Poisson bracket
can be written as a sum of the (trivial) linear Poisson bracket
on $(\Gh/\Gz)^*$ and the map $\omega$, where $\omega:\bigl((x,y),
(x',y')\bigr)\mapsto\eta_{\lambda}(r)\bigl(\beta(x,y'),\beta(x',y)
\bigr)$ is a Lie algebra cocycle on $\Gh/\Gz$ having values in
$C^{\infty}({\Gh^*}/{\Gz^{\bot}})$.  We then obtain the group
cocycle $\sigma$ above, by ``integrating'' $\omega$.  In a more
general setting, this procedure of finding a group cocycle from
a Poisson bracket is discussed in \cite{BJKp1} (See, in particular,
the discussion from Theorem 2.2 through Proposition 3.3 in that
paper.).
\end{rem}

In addition to giving us the group cocycle $\sigma$, the
Poisson bracket strongly suggests us to work with the
$(x,y;r)$ variables, where $(x,y)\in H/Z$ and $r\in
{\Gh^*}/{\Gz^{\bot}}$.  Dual space to $H/Z$ is $(\Gh/\Gz)^*
=\Gz^{\bot}$, whose typical element is denoted by $(p,q)$.
Let us take this suggestion and break the group $G$ into
two, obtaining the following matched pair:

\begin{defn}\label{mpair}
Let $G_1$ and $G_2$ be subgroups of $G$, defined by
$$
G_1=\bigl\{(0,0,r):r\in\mathbb{R}\bigr\},\qquad
G_2=\bigl\{(p,q,0):p,q\in\mathbb{R}^n\bigr\}.
$$
Clearly, as a space $G\cong G_2\times G_1$.  Moreover,
$G_1$ and $G_2$ are closed subgroups of $G$, such that
$G_1\cap G_2=\bigl\{(0,0,0)\bigr\}$.  And, any element
$(p,q,r)$ of $G$ can be (uniquely) expressed as a
product: $(p,q,r)=(0,0,r)(p,q,0)$, with $(0,0,r)\in
G_1$ and $(p,q,0)\in G_2$.  In other words, the groups
$G_1$ and $G_2$ form a {\em matched pair\/}.

From this, we naturally obtain the group actions
$\alpha:G_1\times G_2\to G_2$ and $\gamma:G_2
\times G_1\to G_1$, defined by
$$
\alpha_r(p,q):=(e^{-\lambda r}p,e^{-\lambda r}q),
\qquad\gamma_{(p,q)}(r):=r.
$$
Here we are using the obvious identification of $(p,q)$
with $(p,q,0)$, and similarly for $r$ and $(0,0,r)$.
Note that these actions are defined so that we have: 
$\bigl(\alpha_r(p,q)\bigr)\bigl(\gamma_{(p,q)}(r)\bigr)
=(e^{-\lambda r}p,e^{-\lambda r}q,0)(0,0,r)=(p,q,r)$.
\end{defn}

Let us now convert the informations we obtained so far
into the language of Hilbert space operators and operator
algebras.  To begin with, let us fix a Lebesgue measure
on $H(=\Gh)$, which is the Haar measure for $H$.  And
on $G(=\Gg)$, which is considered as the dual vector
space of $H$, we give the dual Lebesgue measure (As noted
earlier, this will be again the Haar measure for $G$.).
They are chosen so that the Fourier transform becomes
the unitary operator (from $L^2(H)$ to $L^2(G)$),
and the Fourier inversion theorem holds.  Similarly,
``partial'' Fourier transform can be considered: For
instance, between functions in the $(p,q;r)$ variables
and those in the $(x,y;r)$ variables.  See Remark 1.7
of \cite{BJKp2}.

Following Baaj and Skandalis \cite{BS}, the information
about the groups $G_1$ and $G_2$ can be incorporated
into certain multiplicative unitary operators $X$ and $Y$.
The result is given below.  Note that we are also expressing
our operators in the $(x,y;r)$ variables, so that we can
later work within that setting.

\begin{prop}\label{XY}
Let $X\in{\mathcal B}\bigl(L^2(G_1\times G_1)\bigr)$ and
$Y\in{\mathcal B}\bigl(L^2(G_2\times G_2)\bigr)$ be defined
such that for $\xi\in L^2(G_1\times G_1)$ and $\zeta\in 
L^2(G_2\times G_2)$,
$$
X\xi(r;r')=\xi(r+r';r'),\qquad Y\zeta(p,q;p',q')
=\zeta(p+p',q+q';p',q').
$$
They are multiplicative unitary operators.  Meanwhile, by
Fourier transform, $Y$ can be expressed as an operator
contained in ${\mathcal B}\bigl(L^2(H/Z\times H/Z)\bigr)$,
in the $(x,y)$ variables.  This means that we are regarding
${\mathcal F}^{-1}Y{\mathcal F}$ as same as $Y$, for convenience.
It then reads:
$$
Y\zeta(x,y;x',y')=\zeta(x,y;x'-x,y'-y),\quad\zeta\in L^2(H/Z).
$$
We have:
\begin{align}
C_0(G_1)&\cong\overline{\bigl\{(\omega\otimes\operatorname{id})
(X):\omega\in{\mathcal B}(L^2(G_1))_*\bigr\}}^{\|\ \|}\subseteq
{\mathcal B}\bigl(L^2(G_1)\bigr),  \notag \\
C^*(G_1)&\cong\overline{\bigl\{(\operatorname{id}\otimes\omega)
(X):\omega\in{\mathcal B}(L^2(G_1))_*\bigr\}}^{\|\ \|}\subseteq
{\mathcal B}\bigl(L^2(G_1)\bigr),  \notag \\
C_0(G_2)&\cong C^*(H/Z)\cong\overline{\bigl\{(\omega\otimes
\operatorname{id})(Y):\omega\in{\mathcal B}(L^2(H/Z))_*
\bigr\}}^{\|\ \|}\subseteq{\mathcal B}\bigl(L^2(H/Z)\bigr),
\notag  \\
C^*(G_2)&\cong C_0(H/Z)\cong\overline{\bigl\{(\operatorname
{id}\otimes\omega)(Y):\omega\in{\mathcal B}(L^2(H/Z))_*
\bigr\}}^{\|\ \|}\subseteq{\mathcal B}\bigl(L^2(H/Z)\bigr).
\notag
\end{align}
\end{prop}

\begin{proof}
We are just following \cite{BS}.  For the expression of $Y\in
{\mathcal B}\bigl(L^2(H/Z)\bigr)$, we just use the Fourier
inversion theorem.  Since the groups are abelian, all the
computations are quite simple.
\end{proof}

By the result of Proposition \ref{XY}, a function $f\in 
C_0(G_1)$ is considered same as the multiplication operator
$L_f\in{\mathcal B}\bigl(L^2(G_1)\bigr)$, defined by $L_f
\xi(r)=f(r)\xi(r)$.  Similar for $g\in C_0(G_2)$, which is
also considered as the multiplication operator $\lambda_g
\in{\mathcal B}\bigl(L^2(G_2)\bigr)$.  In the $(x,y)$
variables, this is equivalent to saying that for $g\in 
C_c(H/Z)\subseteq C^*(H/Z)$, the operator $\lambda_g
\in{\mathcal B}\bigl(L^2(H/Z)\bigr)$ is such that for
$\zeta\in L^2(H/Z)$, we have: $\lambda_g\zeta(x,y)
=\int g(\tilde{x},\tilde{y})\zeta(x-\tilde{x},
y-\tilde{y})\,d\tilde{x}d\tilde{y}$.

At the level of the $C^*$-algebras $C_0(G_1)$ and $C_0(G_2)$,
the group actions $\alpha$ and $\gamma$ we obtained earlier
(though $\gamma$ is trivial) are expressed as coactions
$\alpha:C_0(G_2)\to M\bigl(C_0(G_2)\otimes C_0(G_1)\bigr)$
and $\gamma:C_0(G_1)\to M\bigl(C_0(G_2)\otimes C_0(G_1)\bigr)$,
given by
\begin{align}
\alpha(g)(p,q;r)&=g(e^{-\lambda r}p,e^{-\lambda r}q)
=g\bigl(\alpha_r(p,q)\bigr),  \notag \\
\gamma(f)(p,q;r)&=f(r)=f\bigl(\gamma_{(p,q)}(r)\bigr).  \notag
\end{align}
Furthermore, the coactions $\alpha$ and $\gamma$ can be realized
using a certain unitary operator $Z$, as follows:

\begin{prop}\label{Z}
Let $Z\in{\mathcal B}\bigl(L^2(G)\bigr)={\mathcal B}\bigl(L^2
(G_2\times G_1)\bigr)$ be defined by
$$
Z\xi(p,q;r)=(e^{-\lambda r})^n\xi(e^{-\lambda r}p,e^{-\lambda r}q;r).
$$
Then we have, for $g\in C_0(G_2)$ and $f\in C_0(G_1)$,
$$
Z(\lambda_g\otimes1)Z^*=(\lambda\otimes L)\bigl(\alpha(g)\bigr),
\qquad Z(1\otimes L_f)Z^*=(\lambda\otimes L)\bigl(\gamma(f)\bigr).
$$
\end{prop}

\begin{proof}
A straightforward computation shows that for $\xi\in L^2(G)$,
$$
Z(\lambda_g\otimes1)Z^*\xi(p,q,r)=g(e^{-\lambda r}p,e^{-\lambda r}q)
\xi(p,q,r)=(\lambda\otimes L)\bigl(\alpha(g)\bigr)\xi(p,q,r).
$$
And similarly, $Z(1\otimes L_f)Z^*=(\lambda\otimes L)\bigl(\gamma(f)
\bigr)$.
\end{proof}

\begin{rem}
By using the operator realizations $g=\lambda_g$ and $f=L_f$,
as well as $\alpha(g)=(\lambda\otimes L)\bigl(\alpha(g)\bigr)$
and $\gamma(f)=(\lambda\otimes L)\bigl(\gamma(f)\bigr)$, we
may simply write the above result as: $Z(g\otimes1)Z^*=\alpha(g)$
and $Z(1\otimes f)Z^*=\gamma(f)$. 
\end{rem}

Since we prefer to work with the $(x,y;r)$ variables, it will be
more convenient to introduce the Hilbert space ${\mathcal H}:=
L^2(H/Z\times G_1)$, consisting of the $L^2$-functions in the
$(x,y;r)$ variables.  Then by considering that $G_2=(H/Z)^*$,
or equivalently that $C_0(G_2)\cong C^*(H/Z)$, we may as well
regard the coactions $\alpha$ and $\gamma$ to be on $C^*(H/Z)$
and $C_0(G_1)$ (In that case, the definitions of $\alpha$ and
$\gamma$ should be modified accordingly.).  In this setting,
the operator $Z$ will become $Z\in{\mathcal B}({\mathcal H})$,
defined by
$$
Z\xi(x,y;r)=(e^{\lambda r})^n\xi(e^{\lambda r}x,e^{\lambda r}y;r).
$$

The multiplicative unitary operator associated to the matched pair
$(G_1,G_2)$ is given in the next proposition.  It gives us two
$C^*$-bialgebras, which are the bicrossed product algebras coming
from the matched pair.  Note that at this moment, the cocycle is not
considered yet.

\begin{prop}\label{V}
Let $V\in{\mathcal B}({\mathcal H}\otimes{\mathcal H})={\mathcal B}
\bigl(L^2(H/Z\times G_1\times H/Z\times G_1)\bigr)$ be the unitary
operator defined by
$$
V=(Z_{12}X_{24}Z^*_{12})Y_{13},
$$
where we are using the standard leg notation.  More specifically,
\begin{align}
&V\xi(x,y,r;x',y',r')   \notag \\
&=(e^{-\lambda r'})^n\xi(e^{-\lambda r'}x,e^{-\lambda r'}y,r+r';
x'-e^{-\lambda r'}x,y'-e^{-\lambda r'}y,r').  \notag
\end{align}

It is multiplicative, and by considering the ``left [and right]
slices'' of $V$, we obtain the following two $C^*$-algebras, as
contained in the operator algebra ${\mathcal B}({\mathcal H})$:
\begin{align}
A_V&=\bigl\{(\omega\otimes\operatorname{id}_{\mathcal H})(V):
\omega\in{\mathcal B}({\mathcal H})_*\bigr\}\cong C_0(G_1)
\rtimes_{\gamma}(H/Z)  \notag \\
\hat{A}_V&=\bigl\{(\operatorname{id}_{\mathcal H}\otimes\omega)
(V):\omega\in{\mathcal B}({\mathcal H})_*\bigr\}\cong C_0(H/Z)
\rtimes_{\alpha}G_1.  \notag
\end{align}
Here the actions are defined by $\gamma_{(x,y)}(r):=r$ (trivial
action), and $\alpha_r(x,y):=(e^{\lambda r}x,e^{\lambda r}y)$.
The algebras $A_V$ and $\hat{A}_V$ are actually $C^*$-bialgebras,
whose comultiplications are given by $\Delta_V(a)=V(a\otimes1)
V^*$ for $a\in A_V$, and $\hat{\Delta}_V(b)=V^*(1\otimes b)V$
for $b\in\hat{A}_V$.
\end{prop}

\begin{rem}
The choice of the operator $V$ is suggested from Section 8 of
\cite{BS}, where discussions are given on obtaining multiplicative
unitary operators from a matched pair ({\em couple assorti\/}).
The point is that the operators $X$ and $Y$ encode the groups
$G_1$ and $G_2$, while the operator $Z$ carries the information
about the actions $\alpha$ and $\gamma$.  Indeed, the statement
above concerning the characterizations for $A_V$ and $\hat{A}_V$
is a fairly general result.
\end{rem}

\begin{proof}
We skip the proof that $V$ is multiplicative, though a direct
verification of the multiplicativity is not really difficult.
Instead, we point out that $V$ is actually a degenerate case
of the multiplicative unitary operator $U$ given in \cite{BJKp2}
(See proof of Proposition \ref{Vtheta} below.).  Once it is
known that $V$ is indeed multiplicative, the general theory
assures us that we have a pair of $C^*$-bialgebras $A_V$ and
$\hat{A}_V$, contained in ${\mathcal B}({\mathcal H})$.

Meanwhile, it is also not difficult to show directly that as
a $C^*$-algebra, we have: $A_V\cong\overline{L\bigl(C_c(H/Z
\times G_1)\bigr)}^{\|\ \|}$, where
$$
L_f\xi(x,y;r)=\int f(\tilde{x},\tilde{y};r)\xi(x-\tilde{x},
y-\tilde{y};r)\,d\tilde{x}d\tilde{y},
$$
for $f\in C_c(H/Z\times G_1)$ and $\xi\in{\mathcal H}$.
It then follows immediately that: $A_V\cong C_0(G_1)\rtimes
(H/Z)$, which is the crossed product algebra with the trivial
action.

Similarly, we can also show that $\hat{A}_V\cong\overline
{\rho\bigl(C_c(H/Z\times G_1)\bigr)}^{\|\ \|}$, where
\begin{align}
\rho_f\xi(x,y;r)&=\int f(x,y;\tilde{r})\xi(e^{\lambda
\tilde{r}}x,e^{\lambda\tilde{r}}y;r-\tilde{r})\,d\tilde{r}
\notag \\
&=\int f(x,y;\tilde{r})\xi\bigl(\alpha_{\tilde{r}}(x,y);
r-\tilde{r})\,d\tilde{r},
\notag
\end{align}
for $f\in C_c(H/Z\times G_1)$ and $\xi\in{\mathcal H}$.
From this, we see that as a $C^*$-algebra, $\hat{A}_V\cong 
C_0(H/Z)\rtimes_{\alpha}G_1$.

We did not provide explicit computations here, but see the
comments made in proof of Proposition \ref{Vtheta}.  The
$C^*$-algebras $A_V$ and $\hat{A}_V$ considered here are
actually degenerate cases of the ones in that proposition.
See also the proof of Proposition \ref{case3}, which has
a similar result.
\end{proof}

Now that we have the multiplicative unitary operator $V$ encoding
the matched pair $(G_1,G_2)$, our final task is to incorporate the
cocycle term.  In the setting of multiplicative unitary operators,
we need to look for a function $\Theta:(H/Z\times G_1)\times(H/Z
\times G_1)\to\mathbb{T}$ such that $V\Theta$ is still multiplicative
(See Section 8 of \cite{BS}.).  Note here that we are regarding
$\Theta$ as a unitary  operator such that $\Theta\xi(x,y,r;
x',y',r')=\Theta(x,y,r;x',y',r')\xi(x,y,r;x',y',r')$.  Motivated
by Proposition \ref{cocycle}, let us take $\Theta$ to be the map,
$$
\Theta(x,y,r;x',y',r'):=\bar{e}\bigl[\eta_{\lambda}(r')
\beta(x,y')\bigr].
$$
As the next proposition shows, we obtain in this way a multiplicative
unitary operator $V_{\Theta}=V\Theta$.  It determines a pair of
$C^*$-bialgebras that are realized as cocycle bicrossed products.

\begin{prop}\label{Vtheta}
[Quantization of Case (2).]
Let $\Theta$ and $V$ be as above.  Then the operator $V_{\Theta}
:=V\Theta\in{\mathcal B}({\mathcal H}\otimes{\mathcal H})$ is
a multiplicative unitary operator.  Specifically,
\begin{align}
V_{\Theta}\xi(x,y,r;x',y',r')&=(e^{-\lambda r'})^n\bar{e}\bigl[
\eta_{\lambda}(r')\beta(e^{-\lambda r'}x,y'-e^{-\lambda r'}y)\bigr]
\notag \\
&\quad\xi(e^{-\lambda r'}x,e^{-\lambda r'}y,r+r';x'-e^{-\lambda r'}x,
y'-e^{-\lambda r'}y,r').  \notag
\end{align}
The $C^*$-bialgebras associated with $V_{\Theta}$ are:
$$
A\cong C_0(G_1)\rtimes_{\gamma}^{\sigma}(H/Z),\quad {\text {and }}
\quad \hat{A}\cong C_0(H/Z)\rtimes_{\alpha}G_1,
$$
together with the comultiplications: $\Delta(a)=V_{\Theta}
(a\otimes1)V_{\Theta}^*$ for $a\in A$, and $\hat{\Delta}(b)
=V_{\Theta}^*(1\otimes b)V_{\Theta}$ for $b\in\hat{A}$.
\end{prop}

\begin{proof}
The operator $V_{\Theta}$ coincides with the multiplicative unitary
operator $U$ obtained in Proposition 3.1 of \cite{BJKp2}.  We will
refer to that paper for the proof of the multiplicativity.  If
$\Theta\equiv1$, the operator $V_{\Theta}$ degenerates into $V$
given in Proposition \ref{V}, giving us the proof of its
multiplicativity we skipped.

As for the characterization of the $C^*$-algebra $A$ as a twisted
crossed product algebra (in the sense of \cite{PR}), see Proposition
2.2 of \cite{BJKppha}, as well as \cite{BJKp2}.  As noted earlier,
$\gamma$ is actually a trivial cocycle, while $\sigma$ is the group
cocycle for $H/Z$ defined in Proposition \ref{cocycle}.  In case
$\sigma\equiv1$ (corresponding to $\Theta\equiv1$), it will
degenerate to $A_V\cong C_0(G_1)\rtimes(H/Z)$ in Proposition
\ref{V}.

The characterization for the $C^*$-algebra $\hat{A}$ can be found
in Proposition~2.2 of \cite{BJKqhg}.  Note that $\hat{A}$ does not
change from the case without the cocycle, given in Proposition
\ref{V} (so $\hat{A}\cong\hat{A}_V$).  Only its comultiplication
changes, by $\hat{\Delta}(b)=\Theta^*\hat{\Delta}_V(b)\Theta$.
\end{proof}

Our approach was different, but since we obtained the same
multiplicative unitary operator as in \cite{BJKp2}, we can
use the result of that paper (as well as \cite{BJKppha}) to
construct the rest of the quantum group structure for $(A,
\Delta)$.  It is a ``quantized $C_0(G)$'', as well as a
``quantized $C^*(H)$'' (Note that $A\cong C^*(H)$, if
$\lambda=0$.).  It is the quantization of the Poisson--Lie
group $G$ given in Proposition \ref{dualLiegroups} (2).

Reflecting the fact that the group $G$ was non-unimodular,
the quantum group $(A,\Delta)$ turns out to be also
non-unimodular (See \cite{BJKppha}.).  See also \cite{BJKhj}
and \cite{BJKppdress}, where we take advantage of the close
relationship between $(G,H)$ and $(A,\Delta)$ to discuss
the representation theory of $A$.  For instance, a
``quasitriangular, quantum $R$-matrix'' type operator
can be found, corresponding to the classical $r$-matrix
given in Appendix (Section 4).

Meanwhile, $(\hat{A},\hat{\Delta})$ is the dual quantum group
of $(A,\Delta)$.  This is studied in \cite{BJKqhg}, and may be
regarded as a ``quantized $C_0(H)$'', or a ``quantized $C^*(G)$''.
As $H$ is unimodular, so is $(\hat{A},\hat{\Delta})$. This is
a ``quantum Heisenberg group'', but is different from the one
constructed in \cite{SZ} or in \cite{VD}.  See below.

\bigskip

Since we are satisfied with Case (2), let us now turn our attention
to Case (1).  Consider the Poisson--Lie group $G$ and
the Poisson bracket on it as described in Proposition \ref
{dualLiegroups} (1).  Noting the similarity with Case (2),
and with a slight modification of the procedure, we obtain
the following:

\begin{prop}\label{case1}
[Quantization of Case (1).]
\begin{enumerate}
\item Let $G_1$ and $G_2$ be defined by
$$
G_1=\{r:r\in\mathbb{R}\},\qquad 
G_2=\bigl\{(p,q):p,q\in\mathbb{R}^n\bigr\}.
$$
Consider also the group actions $\alpha:G_1\times G_2\to G_2$ and
$\gamma:G_2\times G_1\to G_1$, given by
$$
\alpha_r(p,q):=(e^{-\lambda r}p,e^{\lambda r}q),\qquad 
\gamma_{(p,q)}(r):=r.
$$
In this way, we obtain the matched pair $(G_1,G_2)$.
\item Let $X\in{\mathcal B}\bigl(L^2(G_1\times G_1)\bigr)$ and 
$Y\in{\mathcal B}\bigl(L^2(H/Z\times H/Z)\bigr)$ be the operators
defined by $X\xi(r;r')=\xi(r+r';r')$, and by $Y\xi(x,y;x',y')
=\xi(x,y;x'-x,y'-y)$.  In addition, let $Z\in{\mathcal B}
\bigl(L^2(H/Z\times G_1)\bigr)$ be such that $Z\xi(x,y;r)
=\xi(e^{\lambda r}x,e^{-\lambda r}y;r)$, and let $\Theta
(x,y,r;x',y',r'):=\bar{e}\bigl[r'\beta(x,y')\bigr]$, which is
considered as a unitary operator.  Then $V_{\Theta}:=(Z_{12}
X_{24}Z^*_{12}Y_{13})\Theta$ is a multiplicative unitary
operator contained in ${\mathcal B}({\mathcal H}\otimes
{\mathcal H})={\mathcal B}\bigl(L^2(H/Z\times G_1\times H/Z
\times G_1)\bigr)$.
\item The $C^*$-bialgebras associated with $V_{\Theta}$ are:
$$
A\cong C_0(G_1)\rtimes_{\gamma}^{\sigma}(H/Z)\cong C^*(H),
\quad {\text {and }}\quad \hat{A}\cong C_0(H/Z)\rtimes_{\alpha}
G_1,
$$
together with the comultiplications : $\Delta(a):=V_{\Theta}
(a\otimes1)V_{\Theta}^*$ for $a\in A$, and $\hat{\Delta}(b):=
V_{\Theta}^*(1\otimes b)V_{\Theta}$ for $b\in\hat{A}$.
\end{enumerate}
\end{prop}

\begin{rem}
The proof is done in exactly the same way as in the earlier
part of this section, concerning Case (2). Similarly to
Case (2), the $C^*$-algebra $A$ is isomorphic to a twisted
crossed product algebra, with the twisting cocycle $\sigma^r
\bigl((x,y),(x',y')\bigr):=\bar{e}\bigl[r\beta(x,y')\bigr]$.
But by using (partial) Fourier transform, it can be shown
easily that $A\cong C^*(H)$.  This does not hold in Case (2).
It reflects the fact that in Case (1), unlike in Case (2),
the Poisson bracket we began with is a linear Poisson bracket
(dual to the Lie bracket on $H$).
\end{rem}

Again the approach was different, but we point out here
that $(\hat{A},\hat{\Delta})$ from Proposition \ref{case1}
is isomorphic to the example given by Van Daele in \cite{VD},
as well as to the one given by Szymczak and Zakrzewski
\cite{SZ}.  Meanwhile, $(A,\Delta)$ above is its dual
counterpart.  This is actually the special case of the
example considered by Rieffel in \cite{Rf5} (See also
\S3.2 below.).  These three examples are considered
to be among the pioneering works on non-compact quantum
groups.  We will refer to these other papers for the
construction of the rest of quantum group structures for
both $(\hat{A},\hat{\Delta})$ and $(A,\Delta)$.

Representation theory for Case (1) has not been done
in the literatures, but actually, due to the fact that
it corresponds to a linear Poisson bracket and also
to a triangular classical $r$-matrix (see Appendix:
Section 4), it is much simpler than that of Case (2).
Meanwhile, reflecting the fact that both $H$ and $G$ are
unimodular, both $(\hat{A},\hat{\Delta})$ and $(A,\Delta)$
for Case (1) turns out to be unimodular (i.\,e. their
Haar weights are both right and left invariant).

\subsection{Case (3)}
Let us now consider the case of the Lie group $G$ and
the Poisson bracket on it as described in Proposition
\ref{dualLiegroups} (3).

Analogously to Definition \ref{mpair} and Proposition
\ref{case1} (1), we will begin with the matched pair
$(G_1,G_2)$.  Here, the groups are
$$
G_1=\{r:r\in\mathbb{R}\}\quad {\text {and }}\quad 
G_2=\bigl\{(p,q):p,q\in\mathbb{R}^n\bigr\},
$$
together with the group actions $\alpha:G_1\times G_2\to 
G_2$ and $\gamma:G_2\times G_1\to G_1$, given by
$\alpha_r(p,q):=\left(p,q-r\sum_{i,j}J_{ij}p_i\mathbf{q_j}
\right)$ and $\gamma_{(p,q)}(r):=r$.

Note that, as before, $G\cong G_1\times G_2$ as a space,
while $G_1$ and $G_2$ may be regarded as closed subgroups
of $G$ such that $G_1\cap G_2=\bigl\{(0,0,0)\bigr\}$.
This is done by viewing $(0,0,r)$ and $(p,q,0)$ as same
as $r\in G_1$ and $(p,q)\in G_2$, respectively.  Any
element of $G$ can be (uniquely) expressed as a product:
$(p,q,r)=(0,0,r)(p,q,0)$.  The actions are defined so
that we have: $\bigl(\alpha_r(p,q)\bigr)(\gamma_{(p,q)}(r)
\bigr)=(p,q,r)$.

We will again work with the $(x,y;r)$ variables, in $H/Z
\times G_1$.  So the multiplicative unitary operators
associated with the groups $G_1$ and $G_2$ are $X\in
{\mathcal B}\bigl(L^2(G_1\times G_1)\bigr)$ and $Y\in
{\mathcal B}\bigl(L^2(H/Z\times H/Z)\bigr)$, defined by
$X\xi(r;r')=\xi(r+r';r')$, for $\xi\in L^2(G_1\times G_1)$,
and $Y\xi(x,y;x',y')=\xi(x,y;x'-x,y'-y)$, for $\xi\in L^2
(H/Z\times G_1)$.

The operator encoding the group actions $\alpha$ and $\gamma$
is $Z\in{\mathcal B}\bigl(L^2(G_2\times G_1)\bigr)$, defined
by $Z\xi(p,q;r)=\xi\left(p,q-r\sum_{i,j}J_{ij}p_i\mathbf{q_j};
r\right)$.  By using partial Fourier transform and the Fourier
inversion theorem, we see that it is equivalent to the following
(same-named) operator $Z\in{\mathcal B}\bigl(L^2(H/Z\times G_1)
\bigr)$:
$$
Z\xi(x,y;r)=\xi\left(x+r\sum_{i,j}J_{ij}y_j\mathbf{x_i},y;
r\right),\quad {\text {for }}\quad \xi\in L^2(H/Z\times G_1).
$$
All this is again very much similar to Propositions \ref{XY}
and \ref{Z}, as well as Proposition \ref{case1} (2).

Using the same strategy as in Proposition \ref{Vtheta} or
in Proposition \ref{case1}, we define our multiplicative
unitary operator $V_{\Theta}$, as follows.  In particular,
the definition of the cocycle $\Theta$ comes directly from
the expression of the Poisson bracket given in Proposition 
\ref{dualLiegroups} (3).  In the below, ${\mathcal H}$ denotes
the Hilbert space $L^2(H/Z\times G_1)$.

\begin{prop}\label{case3}
[Quantization of Case (3).]
Define the unitary operator $V\in{\mathcal B}({\mathcal H}
\otimes{\mathcal H})={\mathcal B}\bigl(L^2(H/Z\times G_1\times 
H/Z\times G_1)\bigr)$, by $V=(Z_{12}X_{24}Z^*_{12})Y_{13}$.
It is multiplicative.  And let $\Theta(x,y,r;x',y',r'):=
\bar{e}\bigl[r'\beta(x,y')\bigr]\bar{e}\left[\frac{{r'}^2}{2}
\sum_{i,j}J_{ij}y_jy'_i\right]$, considered as a unitary
operator contained in ${\mathcal B}({\mathcal H}\otimes
{\mathcal H})$.

Then the function $\Theta$ is a cocycle for $V$.  In this way,
we obtain a multiplicative unitary operator $V_{\Theta}:=
V\Theta\in{\mathcal B}({\mathcal H}\otimes{\mathcal H})$.
Specifically,
\begin{align}
&V_{\Theta}\xi(x,y,r;x',y',r')  \notag \\
&=e\left[\frac{{r'}^2}{2}
\sum_{i,j}J_{ij}y_j(y'_i-y_i)\right]\bar{e}\bigl[r'\beta
(x,y'-y)\bigr]  \notag \\
&\qquad\xi\left(x-r'\sum_{i,j}J_{ij}y_j\mathbf{x_i},y,r+r';
x'-x+r'\sum_{i,j}J_{ij}y_j\mathbf{x_i},y'-y,r'\right).
\notag
\end{align}
The $C^*$-bialgebras associated with $V_{\Theta}$ are:
$$
S\cong C_0(G_1)\rtimes_{\gamma}^{\sigma}(H/Z),\quad 
{\text {and }}\quad \hat{S}\cong C_0(H/Z)\rtimes_{\alpha} G_1,
$$
together with the comultiplications : $\Delta(a):=V_{\Theta}
(a\otimes1){V_{\Theta}}^*$ for $a\in S$, and $\hat{\Delta}(b):=
{V_{\Theta}}^*(1\otimes b)V_{\Theta}$ for $b\in\hat{S}$.  Here,
$\sigma:r\mapsto\sigma^r$ is a continuous field of cocycles
such that $\sigma^r\bigl((x,y),(x',y')\bigr)=\bar{e}\left[\frac
{{r}^2}{2}\sum_{i,j}J_{ij}y_jy'_i\right]\bar{e}\bigl[r\beta
(x,y')\bigr]$.
\end{prop}

\begin{proof}
The multiplicativity of $V$ is a consequence of the fact that
$(G_1,G_2)$ forms a matched pair, or equivalently, that $G$
is a group.  The function $\Theta$ is a cocycle for $V$, since
$V_{\Theta}$ is also multiplicative.  The verification of
the pentagon equation, $W_{12}W_{13}W_{23}=W_{23}W_{12}$
for $W=V_{\Theta}$, is straightforward.

As usual, the $C^*$-bialgebras associated with $V_{\Theta}$
are obtained by
\begin{align}
S&=\bigl\{(\omega\otimes\operatorname{id}_{\mathcal H})
(V_{\Theta}):\omega\in{\mathcal B}({\mathcal H})_*\bigr\}
\bigl(\subseteq{\mathcal B}({\mathcal H})\bigr),  \notag \\
\hat{S}&=\bigl\{(\operatorname{id}_{\mathcal H}\otimes\omega)
(V_{\Theta}):\omega\in{\mathcal B}({\mathcal H})_*\bigr\}
\bigl(\subseteq{\mathcal B}({\mathcal H})\bigr).  \notag
\end{align}
To see the specific $C^*$-algebra realization of $S$, consider
its typical element $(\omega\otimes\operatorname{id}_{\mathcal H})
(V_{\Theta})$, where $\omega\in{\mathcal B}({\mathcal H})_*$.
Without loss of generality, we may assume that $\omega=
\omega_{\xi,\eta}$, for $\xi,\eta\in{\mathcal H}$ (We may even
assume that $\xi$ and $\eta$ are continuous functions having
compact support.).  It is a rather standard notation, and is
defined by $\omega_{\xi,\eta}(T)=\langle T\xi,\eta\rangle$,
for $T\in{\mathcal B}({\mathcal H})$.  It is known that linear
combinations of the $\omega_{\xi,\eta}$ are (norm) dense in
${\mathcal B}({\mathcal H})_*$.  Now for $\zeta\in{\mathcal H}$,
we have:
\begin{align}
&\bigl((\omega_{\xi,\eta}\otimes\operatorname{id}_{\mathcal H})
(V_{\Theta})\bigr)\zeta(x,y,r)  \notag \\
&=\int\bigl(V_{\Theta}(\xi\otimes
\zeta)\bigr)(\tilde{x},\tilde{y},\tilde{r};x,y,r)\overline{\eta
(\tilde{x},\tilde{y},\tilde{r})}\,d\tilde{x}d\tilde{y}d\tilde{r}
\notag \\
&=\int e\left[\frac{{r}^2}{2}\sum_{i,j}J_{ij}\tilde{y}_j
(y_i-\tilde{y}_i)\right]\bar{e}\bigl[r\beta(\tilde{x},
y-\tilde{y})\bigr]\xi\left(\tilde{x}-r\sum_{i,j}J_{ij}\tilde{y}_j
\mathbf{x_i},\tilde{y},\tilde{r}+r\right)
\notag\\
&\qquad\overline{\eta(\tilde{x},\tilde{y},\tilde{r})}\,\zeta
\left(x-\tilde{x}+r\sum_{i,j}J_{ij}\tilde{y}_j\mathbf{x_i},
y-\tilde{y},r\right)\,d\tilde{x}d\tilde{y}d\tilde{r}
\notag \\
&=\int\bar{e}\left[\frac{{r}^2}{2}\sum_{i,j}J_{ij}\tilde{y}_j
(y_i-\tilde{y}_i)\right]\bar{e}\bigl[r\beta(\tilde{x},
y-\tilde{y})\bigr]\xi(\tilde{x},\tilde{y},\tilde{r}+r)
\notag\\
&\qquad\overline{\eta\left(\tilde{x}+r\sum_{i,j}J_{ij}
\tilde{y}_j\mathbf{x_i},\tilde{y},\tilde{r}\right)}\,\zeta
\left(x-\tilde{x},y-\tilde{y},r\right)\,d\tilde{x}d\tilde{y}
d\tilde{r}  \notag \\
&=\int F(\tilde{x},\tilde{y},r)\sigma^r\bigl((\tilde{x},
\tilde{y}),(x-\tilde{x},y-\tilde{y})\bigr)\zeta(x-\tilde{x},
y-\tilde{y},r)\,d\tilde{x}d\tilde{y},
\notag
\end{align}
where $F(x,y,r)=\int\xi(x,y,\tilde{r}+r)\overline{\eta
\left(x+r\sum_{i,j}J_{ij}y_j\mathbf{x_i},y,\tilde{r}
\right)}\,d\tilde{r}$, which is a continuous function
since $\xi$ and $\eta$ are $L^2$-functions.  And $\sigma^r
\bigl((x,y),(x',y')\bigr)=\bar{e}\left[\frac{{r}^2}{2}
\sum_{i,j}J_{ij}y_jy'_i\right]\bar{e}\bigl[r\beta(x,y')
\bigr]$.  It immediately follows from these observations that:
$$
S\cong\overline{\bigl\{(\omega\otimes\operatorname
{id}_{\mathcal H})(V_{\Theta}):\omega\in{\mathcal B}
({\mathcal H})_*\bigr\}}^{\|\ \|}\cong C_0(G_1)
\rtimes_{\gamma}^{\sigma}(H/Z),
$$
which is the twisted crossed product algebra with (trivial)
action $\gamma$, and whose twisting is given by the cocycle
$\sigma:r\mapsto\sigma^r$.

Similar computation as above (and similar also to the case
of $\hat{A}$ in Proposition \ref{Vtheta} and of $\hat{A}_V$
in Proposition \ref{V}) shows that $\hat{S}\cong C_0(H/Z)
\rtimes_{\alpha} G_1$, which is the crossed product algebra
with action $\alpha$, given by $\alpha_r(x,y)
=\left(x+r\sum_{i,j}J_{ij}y_j\mathbf{x_i},y\right)$.
\end{proof}

Essentially, $(S,\Delta)$ is a ``quantized $C^*(H)$'' or
a ``quantized $C_0(G)$''.  For instance, if $J\equiv0$,
then we have: $S\cong C^*(H)$.  Let us also look at the
comultiplication $\Delta$ of $S$ below, which shows that
it reflects the group multiplication law on $G$.

\begin{prop}
For $\phi\in C_c(G)$, define $L_{\phi}\in{\mathcal B}({\mathcal H})$
be defined by
$$
L_{\phi}\zeta(x,y,r):=\int\phi^{\vee}(\tilde{x},\tilde{y},r)
\sigma^r\bigl((\tilde{x},\tilde{y}),(x-\tilde{x},y-\tilde{y})
\bigr)\zeta(x-\tilde{x},y-\tilde{y},r)\,d\tilde{x}d\tilde{y},
$$
where $\sigma$ is the cocycle as in Proposition \ref{case3},
and $\phi^{\vee}$ denotes the (partial) Fourier transform of
$\phi$.  Namely,
$\phi^{\vee}(x,y,r)=\int\phi(p,q,r)e[p\cdot x+q\cdot y]\,dpdq$.
We know from Proposition \ref{case3} that $S\cong\overline
{L\bigl(C_c(G)\bigr)}^{\|\ \|}$, as a $C^*$-algebra.

The comultiplication, $\Delta$, on $S$ is given by $\Delta(a)
=V_{\Theta}(a\otimes1){V_{\Theta}}^*$ for $a\in S$.  For $\phi
\in C_c(G)$, this becomes: $\Delta(L_{\phi})=(L\otimes L)_{\Delta
(\phi)}$, where $\Delta(\phi)\in C_b(G\times G)$ is the function
defined by
$$
\bigl(\Delta(\phi)\bigr)(p,q,r;p',q',r')=\phi\left(p+p',q+q'
+r'\sum_{i,j}J_{ij}p_i\mathbf{q_j},r+r'\right).
$$
\end{prop}

\begin{proof}
Write $L_{\phi}=\int({\mathcal F}^{-1}\phi)(\tilde{x},\tilde{y},
\tilde{z})L_{\tilde{x},\tilde{y},\tilde{z}}\,d\tilde{x}d\tilde{y}
d\tilde{z}$, where ${\mathcal F}^{-1}\phi\in C_c(H)$ is the (inverse)
Fourier transform of $\phi$.  Then $L_{\tilde{x},\tilde{y},\tilde{z}}
\in{\mathcal B}({\mathcal H})$ is such that
$$
L_{\tilde{x},\tilde{y},\tilde{z}}\xi(x,y,r)=\bar{e}(r\tilde{z})
\sigma^r\bigl((\tilde{x},\tilde{y}),(x-\tilde{x},y-\tilde{y})\bigr)
\xi(x-\tilde{x},y-\tilde{y},r).
$$
Comparing with the definition of $L_{\phi}$ given above, we may
regard $L_{\tilde{x},\tilde{y},\tilde{z}}$ as $L_{\tilde{x},\tilde{y},
\tilde{z}}=L_F$, where the function $F\in C_b(G)$ is such that:
$F(p,q,r)=\bar{e}[p\cdot\tilde{x}+q\cdot\tilde{y}+r\tilde{z}]$.
Actually, $L_{\tilde{x},\tilde{y},\tilde{z}}$ is a multiplier,
contained in $M(S)$.  In a sense, the operators $L_{\tilde{x},
\tilde{y},\tilde{z}}$ for $(\tilde{x},\tilde{y},\tilde{z})\in H$,
form the building blocks for the ``regular representation'' $L$ (or
equivalently, for $C^*$-algebra $S$).

For $\zeta\in{\mathcal H}$, we have:
\begin{align}
&\bigl(\Delta(L_{\tilde{x},\tilde{y},\tilde{z}})\bigr)
\zeta(x,y,r;x',y',r')=V_{\Theta}(L_{\tilde{x},\tilde{y},
\tilde{z}}\otimes1){V_{\Theta}}^*\zeta(x,y,r;x',y',r')
\notag \\
&=\bar{e}\bigl[(r+r')\tilde{z}\bigr]\bar{e}\left[\frac
{r^2}{2}\sum_{i,j}J_{ij}\tilde{y}_j(y_i-\tilde{y}_i)\right]
\bar{e}\bigl[r\beta(\tilde{x},y-\tilde{y})\bigr]  \notag \\
&\quad\bar{e}\left[\frac{{r'}^2}{2}\sum_{i,j}J_{ij}\tilde{y}_j
(y'_i-\tilde{y}_i)\right]\bar{e}\bigl[r'\beta(\tilde{x},
y'-\tilde{y})\bigr]\bar{e}\left[rr'\sum_{i,j}J_{ij}\tilde{y}_j
(y_i-\tilde{y}_i)\right]  \notag \\
&\quad\zeta\left(x-\tilde{x}-r'\sum_{i,j}J_{ij}\tilde{y}_j
\mathbf{x_i},y-\tilde{y},r;x'-\tilde{x},y'-\tilde{y},r'
\right).
\notag
\end{align}

Meanwhile, consider $\Delta(F)\in C_b(G\times G)$, given by
\begin{align}
&\bigl(\Delta(F)\bigr)(p,q,r;p',q',r')  \notag \\
&=\bar{e}\left[(p+p')\cdot\tilde{x}+(q+q')\cdot\tilde{y}
+r'\sum_{i,j}J_{ij}p_i\tilde{y}_j+(r+r')\tilde{z}\right].
\notag
\end{align}
Then by a straightforward computation involving Fourier
inversion theorem, we can see that for $\zeta\in{\mathcal H}$:
$$
(L\otimes L)_{\Delta(F)}\zeta(x,y,r;x',y',r')
=\bigl(\Delta(L_{\tilde{x},\tilde{y},\tilde{z}})\bigr)
\zeta(x,y,r;x',y',r').
$$
In other words, $(L\otimes L)_{\Delta(F)}=\Delta(L_F)$.
Remembering the definitions, it follows easily that
$\Delta(L_{\phi})=(L\otimes L)_{\Delta(\phi)}$ for any
$\phi\in C_c(G)$, where $\Delta(\phi)$ is as defined above.
\end{proof}

\begin{rem}
This proposition shows that for $\phi\in C_c(G)$,
the comultiplication sends it to $\Delta(\phi)\in C_b
(G\times G)$, such that
$$
\bigl(\Delta(\phi)\bigr)(p,q,r;p',q',r')=\phi\bigl(
(p,q,r)(p',q',r')\bigr),
$$
preserving the group multiplication law on $G$ as given
in Proposition~\ref{dualLiegroups} (3).  This result
supports our assertion made earlier that $(S,\Delta)$
is a ``quantized $C_0(G)$''.
\end{rem}

At this moment, the $C^*$-bialgebra $(S,\Delta)$ is just
a quantum semi-group.  For it to be properly considered
as a locally compact quantum group, we need further
discussions on maps like antipode or Haar weight (See 
\cite{KuVa} for general theory.).  For this, we may
follow the methods we used earlier in \cite{BJKppha}
or \cite{BJKqhg}, taking advantage of the fact that
$(S,\Delta)$ is a ``quantized $C_0(G)$''.  Meanwhile,
by introducing a deformation parameter, we could
also show that $(S,\Delta)$ is indeed a deformation
quantization of the Poisson--Lie group $G$, in the
direction of its Poisson bracket given in Proposition
\ref{dualLiegroups} (3) [For Case (2), the deformation
quantization is carried out in \cite{BJKp2}.].

In the current paper, though, we will be content
to have carried out our program, and shown a
constructive method of obtaining an appropriate
multiplicative unitary operator for $(S,\Delta)$.

Meanwhile, notice the similarity between our example
$(S,\Delta)$ above and the one constructed by Enock
and Vainerman in Section 6 of \cite{EV}.  The methods
of construction are rather different between the two.
However, looking at the comultiplications and the
cocycles involved, we see a strong resemblance.  What
this means is that the ingredients at the classical
level (informations about the groups $H$ and $G$) are
more or less the same.

On the other hand, there is a very significant difference.
Namely, the example of \cite{EV} has the underlying von 
Neumann algebra isomorphic to the group von Neumann algebra
${\mathcal L}(H)=C^*(H)''$ of $H$.  While in our case,
$S$ is isomorphic to a twisted crossed product algebra:
Unless $J\equiv0$, the $C^*$-algebra $S$ is not isomorphic
to $C^*(H)$.

In the author's opinion, the example $(S,\Delta)$ given
here has a little more merit, considering that its
Poisson--Lie group counterpart and its multiplicative
unitary operator have all been obtained; the relationship
between the Poisson bracket and the cocycle bicrossed
product construction of the multiplicative unitary operator
have been manifested; as well as that the underlying
$C^*$-algebra is built on the framework of twisted crossed
product algebras (more general than ordinary group
$C^*$-algebras or group von Neumann algebras).

\section{Other examples}

In this section, we give more constructions of several
other examples of quantum (semi-)groups.  Just as in the
previous section, each of these examples will have a twisted
crossed product as its underlying $C^*$-algebra, and the
construction can be carried out within the framework of
cocycle bicrossed products.  These examples are actually
slight generalizations of the basic examples given in
Section 2, and they are considered as coming from
Heisenberg-type Lie bialgebras.

Since the new examples will not be fundamentally different
from the basic examples covered in Section 2, we will try
to make discussions rather brief.  We do not even plan to
say much about the Poisson--Lie group counterparts to the
examples.  As before, we will give, for the purpose of
efficiency, just the appropriate multiplicative unitary
operators.  This can be done by giving some small
modifications to the twisting cocycles we had above.
As in Section 2, discussions about correctly establishing
these $C^*$-bialgebras as ($C^*$-algebraic) locally compact
quantum groups will be skipped.

\subsection{``Mixed'' case of Cases (1) and (2).}
Consider the Lie group $G$ defined by the multiplication:
$$
(p,q,r)(p',q',r')=(e^{\lambda r'}p+p',e^{\nu r'}q+q',r+r'),
$$
where $\lambda,\nu\in\mathbb{R}$.  Note that if $\nu=-\lambda$
or $\nu=\lambda$, it coincides with the group $G$ given in
Case (1) or Case (2) of Proposition \ref{dualLiegroups},
respectively.  In fact, the Lie group $G$ above is obtained
as a dual Poisson--Lie group of $(H,\delta_4)$, where
$\delta_4:\Gh\to\Gh\wedge\Gh$ is the cobracket defined by
$\delta_4=\left(\frac{\lambda-\nu}{2\lambda}\right)\delta_1
+\left(\frac{\lambda+\nu}{2\lambda}\right)\delta_2$ [Recall
Definition \ref{PoissonH}.].  In this sense, it is a ``mixed''
case of Case (1) and Case (2) earlier.

To find the quantum counterpart of $G$, or equivalently,
the multiplicative unitary operator for the quantum
(semi-)group, it really boils down to ``changing of the
cocycles''.  So as before, let ${\mathcal H}$ be the Hilbert
space consisting of $L^2$-functions in the $(x,y,r)$ variables.
Also let $\eta_{(\lambda,\nu)}:=\frac{e^{(\lambda+\nu)r}-1}
{\lambda+\nu}$.  Then define $V_{\Theta}\in{\mathcal B}
({\mathcal H}\otimes{\mathcal H})$, given by
\begin{align}
V_{\Theta}\xi(x,y,r;x',y',r')
&=\left(e^{-\frac{(\lambda+\nu)r'}{2}}\right)^n\bar{e}
\bigl[\eta_{(\lambda,\nu)}(r')\beta(e^{-\lambda r'}x,
y'-e^{-\nu r'}y)\bigr]  \notag \\
&\quad\xi(e^{-\lambda r'}x,e^{-\nu r'}y,r+r';
x'-e^{-\lambda r'}x,y'-e^{-\nu r'}y,r').  \notag 
\end{align}
It is obtained following pretty much the same procedure as
in the cases considered in Section 2.  As we see below,
it determines a twisted crossed product algebra whose
twisting cocycle is given by $\sigma^r\bigl((x,y),(x',y')
\bigr):=\bar{e}\bigl[\eta_{(\lambda,\nu)}(r)\beta(x,y')\bigr]$.

\begin{prop}\label{mixed}
Let $V_{\Theta}$ be as in the previous paragraph.  It is
a multiplicative unitary operator.  The $C^*$-bialgebras
associated with $V_{\Theta}$ are:
$$
A\cong C_0(G_1)\rtimes_{\gamma}^{\sigma}(H/Z),\quad 
{\text {and }}\quad \hat{A}\cong C_0(H/Z)\rtimes_{\alpha} G_1,
$$
together with the comultiplications : $\Delta(a):=V_{\Theta}
(a\otimes1){V_{\Theta}}^*$ for $a\in A$, and $\hat{\Delta}(b):=
{V_{\Theta}}^*(1\otimes b)V_{\Theta}$ for $b\in\hat{A}$.  Here,
$G_1$ is the abelian group $G_1=\{r:r\in\mathbb{R}\}$; the
action $\gamma$ is trivial; and the action $\alpha$ is defined
by $\alpha_r(x,y)=(e^{\lambda r}x,e^{\nu r}y)$.  Finally,
$\sigma:r\mapsto\sigma^r$ is a continuous field of cocycles
such that
$$
\sigma^r:H/Z\times H/Z\ni\bigl((x,y),(x',y')\bigr)\mapsto
\bar{e}\bigl[\eta_{(\lambda,\nu)}(r)\beta(x,y')\bigr]
\in\mathbb{T}.
$$
\end{prop}

\begin{proof}
Checking the multiplicativity of $V_{\Theta}$ is straightforward.
To see the realizations of the $C^*$-algebras $A$ and $\hat{A}$,
we use the same method as in the proof of Proposition \ref{case3},
investigating the operators $(\omega\otimes\operatorname{id})
(V_{\Theta})$, $\omega\in{\mathcal B}({\mathcal H})_*$ and
$(\operatorname{id}\otimes\omega')(V_{\Theta})$, $\omega'
\in{\mathcal B}({\mathcal H})_*$.

In this way, we can show that $A\cong\overline{L\bigl(C_c
(H/Z\times G_1)\bigr)}^{\|\ \|}\bigl(\subseteq{\mathcal B}
({\mathcal H})\bigr)$, where $L$ is the ``regular''
representation defined by
$$
L_f\xi(x,y,r)=\int f(\tilde{x},\tilde{y},r)\bar{e}\bigl[
\eta_{(\lambda,\nu)}(r)\beta(\tilde{x},y-\tilde{y})\bigr]
\xi(x-\tilde{x},y-\tilde{y},r)\,d\tilde{x}d\tilde{y}.
$$
Here $f\in C_c(H/Z\times G_1)$ and $\xi\in{\mathcal H}$.
We can see from this observation that the $C^*$-algebra
$A$ is a twisted crossed product algebra, with trivial
action and the twisting cocycle given by $\sigma:r\mapsto
\sigma^r$.

Similarly, we can also show that $\hat{A}\cong\overline
{\rho\bigl(C_c(H/Z\times G_1)\bigr)}^{\|\ \|}\bigl
(\subseteq{\mathcal B}({\mathcal H})\bigr)$, where
$\rho$ is also the regular representation defined by
$$
\rho_f\xi(x,y,r)=\int f(x,y,\tilde{r})\xi(e^{\lambda
\tilde{r}}x,e^{\nu\tilde{r}}y,r-\tilde{r})\,d\tilde{r},
$$
for $f\in C_c(H/Z\times G_1)$ and $\xi\in{\mathcal H}$.
From this, it follows easily that $\hat{A}\cong C_0(H/Z)
\rtimes_{\alpha}G_1$, which is the crossed product
algebra with action $\alpha$.
\end{proof}

Similarly as before, $(A,\Delta)$ is considered as a
``quantized $C^*(H)$'' or a ``quantized $C_0(G)$''.
Further discussion about this case will parallel that
of Case (2).

\subsection{Example of Rieffel's (\cite{Rf5}).}
Let us now allow our group $H$ to have a higher dimensional
center, $Z=\{(0,0,z):z\in\mathbb{R}^m\}$.  Then $H$ will be
now $(2n+m)$-dimensional.  For convenience, let us keep the
same notation and  express the group law on $H$ as
$$
(x,y,z)(x',y',z')=(x+x',y+y',z+z'+\beta(x,y')).
$$
The differences from the definition of $H$ given in Section 1
are that $z$ and $z'$ are now regarded as vectors (for
instance, $z=z_1\mathbf{z_1}+\dots+z_m\mathbf{z_m}$), and
that $\beta(\ ,\ )$ is no longer the inner product.  It will
be understood as a $Z$-valued bilinear map.  The new group
$H$ is still a two-step nilpotent Lie group which closely
resembles the Heisenberg Lie group.  This is actually the
group considered by Rieffel in \cite{Rf5}.

Let $G$ be defined by the multiplication law:
$$
(p,q,r)(p',q',r')=\bigl(\pi(r')p+p',\rho(r')q+q',r+r'\bigr),
$$ 
where $\pi$ and $\rho$ are representations of the group
$G_1=\{(0,0,r):r\in\mathbb{R}^m\}$ on the spaces $\{(p,0,0):
p\in\mathbb{R}^n\}$ and $\{(0,q,0):q\in\mathbb {R}^n\}$,
respectively.  Let us impose the following ``compatibility
condition'' between $\pi,\rho$, and $\beta$, as given by
Rieffel:
\begin{quote}
(Compatibility condition \cite{Rf5}): Assume that $\beta
\bigl(\pi(r)^t x,\rho(r)^t y\bigr)=\beta(x,y)$ and that
$\bigl(\operatorname{det}(\pi(r)\bigr)\bigl(\operatorname
{det}(\rho(r)\bigr)=1$, for all $r,x,y$.
\end{quote}

This compatibility condition makes the group $G$ closely
analogous to our Case (1) earlier, in the sense that $G$
(together with the linear Poisson bracket dual to the Lie
bracket on $H$) becomes the dual Poisson--Lie group of $H$.
If $Z$ is 1-dimensional, the situation will be exactly same
as in Case (1).  Because of this, the quantization can be
carried out in essentially the same way as in Case (1).

We thus obtain the following unitary operator $V_{\Theta}
\in{\mathcal B}({\mathcal H}\otimes{\mathcal H})$, where
${\mathcal H}$ is the Hilbert space consisting of
$L^2$-functions in the $(x,y,r)$ variables:
\begin{align}
&V_{\Theta}\xi(x,y,r;x',y',r')  \notag \\
&=\bar{e}\bigl[r'\cdot\beta(\pi(-r')^t x,y'-\rho(-r')^t y)
\bigr]   \notag \\
&\quad\xi\bigl(\pi(-r')^t x,\rho(-r')^t y,r+r';
x'-\pi(-r')^t x,y'-\rho(-r')^t y,r'\bigr).  \notag 
\end{align}

As before, $V_{\Theta}$ is easily proved to be multiplicative,
and it again determines two $C^*$-bialgebras $(A,\Delta)$ and
$(\hat{A},\hat{\Delta})$ [Result is analogous to Proposition~\ref
{case1}].  As $C^*$-algebras, we will have: $A\cong  C^*(H)$, and
$\hat{A}\cong C_0(H/Z)\rtimes_{\alpha}G_1$, where $\alpha_r(x,y)
=\bigl(\pi(r)^t x,\rho(r)^t y)$.  The $C^*$-algebra $A$ being
isomorphic to the group $C^*$-algebra $C^*(H)$ again reflects
the point that the Poisson bracket on $G$ is linear.

The method was different, but $(A,\Delta)$ and $V_{\Theta}$
obtained in this way are exactly the example constructed by
Rieffel in \cite{Rf5}.  It was really among the first examples
of quantum groups given by deformation quantization process,
and therefore, was the guiding example of all the examples
considered in this work and many others.

\subsection{A two-step solvable Lie group: Non-unimodular case.}
Let $H$ and $G$ be $(2n+m)$-dimensional groups, defined by
the same multiplication laws as in \S3.2.  But this time,
we will no longer require the ``compatibility condition''.
To distinguish the current case from the previous example,
let us assume that $\beta\bigl(\pi(r)^t x,\rho(r)^t y\bigr)
\ne\beta(x,y)$ and that $\bigl(\operatorname{det}(\pi(r)\bigr)
\bigl(\operatorname{det}(\rho(r)\bigr)\ne1$.  So the group $G$
is a non-unimodular, (two-step) solvable Lie group.

Then the setting becomes similar to the example given in
\S3.1.  Therefore, what we need now is to find the cocycle
expression corresponding to $\sigma^r:\bigl((x,y),(x',y')\bigr)
\mapsto\bar{e}\bigl[\eta_{(\lambda,\nu)}(r)\beta(x,y')\bigr]$
of Proposition \ref{mixed}.

Note however that since $Z$ and $G_1$ are higher than
1-dimensional, the counterparts to $\lambda$ and $\nu$
are no longer scalars.  So it is somewhat difficult to
make sense of the expression $\eta_{(\lambda,\nu)}(r)
=\frac{e^{(\lambda+\nu)r}-1}{\lambda+\nu}$.  On the other
hand, we can get around this problem if we only consider
the numerator part of $\eta_{(\lambda,\nu)}(r)$.  This
means that we are changing the Poisson bracket on $G$ by
the factor of $(\lambda+\nu)$.  Since $\lambda$ and $\nu$
are fixed, this modification will not affect the Poisson
duality between $H$ and $G$  (although we do not give any
explicit description of the Poisson bracket here).

\begin{rem}
One drawback is that in so doing, $\lambda$ and $\nu$ lose
some of their flavors as deformation parameters for the
cobracket.  Furthermore, if $\lambda=\nu=0$, we will no
longer have the linear Poisson bracket as before (We instead
obtain a trivial Poisson bracket.).  Nevertheless, since
$\lambda$ and $\nu$ are fixed, non-zero constants (not
considered as parameters), we do not have to worry about
these problems here.
\end{rem}

So let us look for a counterpart to the cocycle $\bar{e}
\bigl[(e^{(\lambda+\nu)r}-1)\beta(x,y')\bigr]$.  For this,
we will take $\bar{e}\bigl[\Sigma\bigl(\beta(\pi(r)^t x,
\rho(r)^t y')-\beta(x,y')\bigr)\bigr]$, where $\Sigma(\ )
:Z\to\mathbb{R}$ is defined by 
$$
\Sigma(z_1\mathbf{z_1}+z_2\mathbf{z_2}+\dots+z_m\mathbf{z_m})
=z_1+z_2+\dots+z_m.
$$  
Note that when $Z$ is 1-dimensional, both cocycles clearly agree.
Using this, let us write down the following unitary operator
$V_{\Theta}\in{\mathcal B}({\mathcal H}\otimes{\mathcal H})$
(Notice the similarity with the one obtained in \S3.1.):
\begin{align}
&V_{\Theta}\xi(x,y,r;x',y',r')   \notag \\
&=\bigl|(\operatorname{det}(\pi(-r'))(\operatorname{det}
(\rho(-r')))\bigr|^{1/2}   \notag \\
&\quad\bar{e}\bigl[\Sigma\bigl(\beta(x,\rho(r')^t y'-y)-\beta
(\pi(-r')^t x,y'-\rho(-r')^t y)\bigr)\bigr]  \notag \\
&\quad\xi\bigl(\pi(-r')^t x,\rho(-r')^t y,r+r';
x'-\pi(-r')^t x, y'-\rho(-r')^t y,r'\bigr).
\notag
\end{align}

It is again multiplicative, and it thus determines
a pair of $C^*$-bialgebras $(A,\Delta)$ and $(\hat{A},
\hat{\Delta})$.  As $C^*$-algebras, we have: $A\cong 
C_0(G_1)\rtimes_{\gamma}^{\sigma}(H/Z)$, and $\hat{A}
\cong C_0(H/Z)\rtimes_{\alpha}G_1$, where $\gamma$ is
the trivial action, $\alpha_r(x,y)=\bigl(\pi(r)^t x,
\rho(r)^t y\bigr)$, and $\sigma:r\mapsto\sigma^r$ is
the continuous field of cocycles given by $\sigma^r
\bigl((x,y),(x',y')\bigr)=\bar{e}\bigl[\Sigma\bigl
(\beta(\pi(r)^t x,\rho(r)^t y')-\beta(x,y')\bigr)
\bigr]$.  All these computations are done following
the same method we have been using so far.

\begin{rem}
Note that the map $\Sigma$ defined above is just the inner
product: $\Sigma(z)=z\cdot\mathbf{1}$, where $\mathbf{1}
=1\mathbf{z_1}+\dots+1\mathbf{z_m}$.  Of course, there is
no particular reason for choosing the vector $\mathbf{1}$
here, and any fixed vector in $Z$ will be sufficient for
our purposes.  Still, they will all give rise to essentially
the same quantum group since we can vary the bilinear map
$\beta(\ ,\ )$ to accommodate the changes. A more significant
observation is that the map $\bigl((x,y),(x',y')\bigr)\mapsto
\beta(\pi(r)^t x,\rho(r)^t y')-\beta(x,y')$ is already an 
additive cocycle having values in $Z$.
\end{rem}

The examples $(A,\Delta)$ and $(\hat{A},\hat{\Delta})$ given
in this subsection are not necessarily very complicated ones.
However, as far as the author knows, these examples have not
been studied before.  On the other hand, they are really
natural generalizations of the examples explored by the author
in his previous papers.  Meanwhile, even though we are not
explicitly investigating Haar weight and the rest of the
quantum group structure maps here, we note that unlike the
example in \S3.2 (or \cite{Rf5}), the quantum group $(A,\Delta)$
given here will be non-unimodular.

\section{Appendix: The classical $r$-matrices and the Poisson
structures on $H$}

In many cases, compatible Poisson brackets on a Lie group
(or equivalently, Lie bialgebra structures) are known to arise
from certain solutions of the classical Yang--Baxter equation
(CYBE), called the {\em classical $r$-matrices\/}.  In this
Appendix, we will first give a very brief background discussion
on classical $r$-matrices (See \cite{Dr}, \cite{CP} for more.).
We will then show that our three basic Lie bialgebra structures
in Definition \ref{PoissonH} are indeed obtained from some
specific classical $r$-matrices.

In general, let $\Gg$ be a Lie algebra and let $r\in\Gg
\otimes\Gg$ be an arbitrary element.  Define a map $\delta_r:
\Gg\to\Gg\otimes\Gg$, by
\begin{equation}\label{(rcoboundary)}
\delta_r(X)=\operatorname{ad}_X(r),\qquad X\in\Gg.
\end{equation}
Then $\delta_r$ is a 1-cocycle on $\Gg$ with values in $\Gg
\otimes\Gg$.  (Actually, it is a coboundary on $\Gg$.).
The following result holds.

\begin{prop} (See \cite{Dr}.)
Let $\Gg$ be a Lie algebra and let $r\in\Gg\otimes\Gg$.  The map 
$\delta_r$ given by equation \eqref{(rcoboundary)} defines a Lie
bialgebra structure on $\Gg$, if and only if the following two
conditions are satisfied:
\begin{itemize}
\item $r^{12}+r^{21}$ is a $\Gg$-invariant element of $\Gg\otimes\Gg$.
\item $[r^{12},r^{13}]+[r^{12},r^{23}]+[r^{13},r^{23}]$ is a 
$\Gg$-invariant element of $\Gg\otimes\Gg\otimes\Gg$.
\end{itemize}
In this case, $\Gg$ is said to be a {\em coboundary Lie bialgebra}.
\end{prop}

The simplest way to satisfy the second condition of the proposition
is to assume that:
$$
[r^{12},r^{13}]+[r^{12},r^{23}]+[r^{13},r^{23}]=0.
$$
This is called the {\em classical Yang--Baxter equation\/} (CYBE).
A solution of the CYBE is called a ``classical $r$-matrix''.  A
coboundary Lie bialgebra structure coming from a solution of the
CYBE is said to be {\em quasitriangular\/}.  If the classical
$r$-matrix further satisfies $r^{12}+r^{21}=0$ (i.\,e.\ it is a
skew solution of the CYBE), it is said to be {\em triangular\/}.
This terminology is closely related with the ``quantum'' situation
and the so-called {\em (quasitriangular/triangular) quantum
universal $R$-matrices\/}.  

\begin{rem}
The quantization problem of triangular and quasitriangular Lie
bialgebras is an important topic in the quantum group theory,
mostly (but not exclusively) at the quantized universal enveloping
algebra (QUE algebra) setting.  Practically, these are the Lie
bialgebras that are more or less expected to be quantized.
Moreover, the ``quantum $R$-matrices'' often play interesting
roles in the representation theory of the quantum group
counterparts to the Poisson--Lie groups (Lie bialgebras).
See \cite{Dr}, \cite{CP}.
\end{rem}

Let us now turn our attention to the $(2n+1)$-dimensional Heisenberg
Lie group $H$ and the Heisenberg Lie algebra $\Gh$, as given in
Section 1.  Consider also the following ``extended'' Heisenberg
Lie algebra:
\begin{defn}
Let $\tilde{\Gh}$ be the $(2n+2)$-dimensional Lie algebra spanned
by $\mathbf{x_i},\mathbf{y_i}(i=1,\dots,n),\mathbf{z},\mathbf{d}$,
with the brackets
$$
[\mathbf{x_i},\mathbf{y_j}]=\delta_{ij}\mathbf{z}, \quad
[\mathbf{d},\mathbf{x_i}]=\mathbf{x_i}, \quad
[\mathbf{d},\mathbf{y_i}]=-\mathbf{y_i}, \quad
\mathbf{z} {\text { is central}}.
$$
The Lie group corresponding to $\tilde{\Gh}$ is the ``extended''
Heisenberg Lie group $\tilde{H}$.  For group multiplication
law on $\tilde{H}$, see Example 3.6 of \cite{BJKp2} or \S2.1
of \cite{BJKhj}.
\end{defn}

For this extended Heisenberg Lie algebra, we can find the following 
solutions of the classical Yang--Baxter equation (CYBE), $r\in\tilde
{\Gh}\otimes\tilde{\Gh}$.  The proofs are straightforward.  It follows
that we now have the Lie bialgebra structures on $\tilde{\Gh}$.

\begin{prop}
\begin{enumerate}
\item Let $r=\lambda(\mathbf{z}\otimes\mathbf{d}-\mathbf{d}\otimes
\mathbf{z})$, $\lambda\ne0$.  Since $\operatorname{span}(\mathbf{z},
\mathbf{d})$ is an abelian subalgebra of $\tilde{\Gh}$, it is easy
to see that $r$ satisfies the CYBE: $[r^{12},r^{13}]+[r^{12},r^{23}]
+[r^{13},r^{23}]=0$.  Since $r$ is also a skew solution, we thus
obtain on $\tilde{\Gh}$ a ``triangular'' Lie bialgebra structure
$\tilde{\delta}_1$, given by $\tilde{\delta}_1(X):=\operatorname
{ad}_X(r)$, $X\in\tilde{\Gh}$.
\begin{align}
\tilde{\delta}_1(\mathbf{x_i})&=\lambda(\mathbf{x_i}\otimes\mathbf{z}
-\mathbf{z}\otimes\mathbf{x_i})=\lambda\mathbf{x_i}\wedge\mathbf{z},
\notag \\
\tilde{\delta}_1(\mathbf{y_i})&=-\lambda\mathbf{y_i}\wedge\mathbf{z},
\quad \tilde{\delta}_1(\mathbf{z})=0,\quad \tilde{\delta}_1(\mathbf{d})=0.
\notag
\end{align}
\item Let $r=2\lambda\bigl(\sum_{i=1}^n(\mathbf{x_i}\otimes
\mathbf{y_i})+\frac{1}{2}(\mathbf{z}\otimes\mathbf{d}+\mathbf{d}
\otimes\mathbf{z})\bigr)$, $\lambda\ne0$.  We can show that $r$
satisfies the CYBE, and also that $r^{12}+r^{21}$ is $\Gh$-invariant.
So we obtain a ``quasitriangular'' Lie bialgebra structure on
$\tilde{\Gh}$ given by the following $\tilde{\delta}_2$:
\begin{align}
\tilde{\delta}_2(\mathbf{x_i})&=\operatorname{ad}_{\mathbf{x_i}}(r)=
\lambda\mathbf{x_i}\wedge\mathbf{z},\quad
\tilde{\delta}_2(\mathbf{y_i})=\operatorname{ad}_{\mathbf{y_i}}(r)=
\lambda\mathbf{y_i}\wedge\mathbf{z}, \notag \\
\tilde{\delta}_2(\mathbf{z})&=\operatorname{ad}_{\mathbf{z}}(r)=0,
\quad \tilde{\delta}_2(\mathbf{d})=\operatorname{ad}_{\mathbf{d}}(r)=0.  
\notag
\end{align}
\end{enumerate}
\end{prop}

Note that $\Gh(\subseteq\tilde{\Gh})$ is a Lie subalgebra of $\tilde{\Gh}$,
and also that $\delta_1$ and $\delta_2$ given in Definition \ref{PoissonH}
are obtained by restricting $\tilde{\delta}_1$ and $\tilde{\delta}_2$ above.
In other words, $(\Gh,\delta_i)$ ($i=1,2$) is a {\em sub-bialgebra\/} of
$(\tilde{\Gh},\tilde{\delta}_i)$ ($i=1,2$), and hence a Lie bialgebra itself.
In this way, we recover the Poisson brackets of Case (1) and Case (2).

For Case (3), see the following proposition (proof is again straightforward).
In this case, we do not need to introduce the extended Heisenberg Lie algebra.
As in Section 1, let $(J_{ij})$ be a skew, $n\times n$ matrix ($n\ge2$).

\begin{prop}
Let $r\in{\Gh}\otimes{\Gh}$ be defined by $r=\sum_{i,j=1}^n J_{ij}\mathbf{x_i}
\otimes\mathbf{x_j}$.  Since $\operatorname{span}(\mathbf{x_i}:i=1,2,\dots,n)
\subseteq\Gh$ is an abelian subalgebra, $r$ clearly satisfies the CYBE.  It
is also a skew solution.  Therefore, we obtain a ``triangular'' Lie bialgebra
structure $\delta_3$ on $\Gh$:
$$
{\delta}_3(\mathbf{x_j})=0, \quad{\delta}_3(\mathbf{y_j})
=\sum_{i=1}^n J_{ij}\mathbf{x_i}\wedge\mathbf{z}, \quad
{\delta}_3(\mathbf{z})=0.  
$$
\end{prop}

Having the knowledge that our Poisson brackets come from certain classical
$r$-matrices is quite useful.  For instance, in \cite{BJKp2} (concerning
Case (2)), we could find an operator $R$, which can be considered as a
``quantum universal $R$-matrix'' (the quantum counterpart to the classical
$r$-matrix).  Using the operator $R$ in \cite{BJKhj}, we could show an
interesting (genuinely quantum) property of ``quasitriangularity'' in the
representation theory of $(A,\Delta)$.


\bibliography{ref}

\bibliographystyle{amsplain}

\end{document}